\documentclass[11pt,reqno]{amsart}
\usepackage{amsmath}
\usepackage{amsfonts}
\usepackage{amssymb,latexsym}
\usepackage{cite}
\usepackage[height=190mm,width=130mm]{geometry}

\theoremstyle{plain}

\newtheorem{proposition}{Proposition}

\theoremstyle{definition}
\newtheorem{definition}{Definition}
\theoremstyle{remark}

\numberwithin{equation}{section}
\input{tcilatex}

\begin{document}
\title[Lucas-Leonardo p-quaternions]{On Companion sequences associated with
Leonardo quaternions: Applications over finite fields}
\author{Diana SAVIN}
\address{Department of Mathematics and Computer Science, Transilvania
University of Brasov, 500091, Romania}
\email{diana.savin@unitbv.ro}
\author{Elif TAN}
\address{Department of Mathematics, Faculty of Science, Ankara University
06100 Tandogan Ankara, Turkey}
\email{etan@ankara.edu.tr}

\begin{abstract}
It is known that the quaternion algebras are central simple algebras and
also clifford algebras. In this paper, we introduce a new class of
quaternions called Lucas-Leonardo p-quaternions and derive several
fundamental properties of these numbers. Furthermore, we investigate some
applications related to companion sequences associated with Leonardo
quaternions. In particular, we determine Lucas-Leonardo quaternions and
Francois quaternions, which are zero divisors and invertible elements in the
quaternion algebra over certain finite fields.
\end{abstract}

\maketitle

\section{Introduction}

Fibonacci sequences and their extensions find numerous applications across
various domains, including arts and sciences. Notably, Fibonacci $p$%
-numbers, a broader variant of Fibonacci numbers, introduced by Stakhov and
Rozin \cite{stakhov} and have found applications in coding theory,
particularly in designing error-correcting codes and data compression
algorithms. For any given integer $p>0,$ the Fibonacci $p$-sequence is
defined by the following recurrence relation%
\begin{equation*}
F_{p,n}=F_{p,n-1}+F_{p,n-p-1},\mbox{   }n>p
\end{equation*}%
with initial values $F_{p,0}=0,$ $F_{p,k}=1$ for $k=1,2,\ldots ,p.$ Its
companion sequence, the Lucas $p$-sequence $\left \{ L_{p,n}\right \}
_{n=0}^{\infty },$ \cite{stakhov} also satisfy the same recurrence relation
but begins with initial values $L_{p,0}=p+1,$ $L_{p,k}=1$ for $k=1,2,\ldots
,p.$ It is clear to see that when $p=1$, the Fibonacci $p$-sequence and the
Lucas $p$-sequence reduce to the classical Fibonacci sequence $\left \{
F_{n}\right \} _{n=0}^{\infty }$ and Lucas sequence $\left \{ L_{n}\right \}
_{n=0}^{\infty }$, respectively. For further information regarding the
Fibonacci $p$-numbers and general second-order linear recurrences, we refer
to \cite{abbad, kocer, stakhov, tuglu}.

There are several non-homogenous extensions of Fibonacci recurrence
relation. One among them is the Leonardo $p$-sequence $\{ \mathcal{L}%
_{p,n}\}_{n=0}^{\infty }$ which was introduced by Tan and Leung \cite%
{tan-leung} and defined by the following non-homogenous relation:
\begin{equation}
\mathcal{L}_{p,n}=\mathcal{L}_{p,n-1}+\mathcal{L}_{p,n-p-1}+p,\mbox{   }n>p
\label{leo}
\end{equation}%
with initial values $\mathcal{L}_{p,0}=\mathcal{L}_{p,1}=\cdots =\mathcal{L}%
_{p,p}=1.$ It is clear to see that when $p=1$, it reduces to the classical
Leonardo sequence$\{ \mathcal{L}_{n}\}_{n=0}^{\infty }.$ The number of
vertices in the $n$-th Leonardo tree is counted by Leonardo numbers. The
Leonardo sequence finds applications in various fields, including
mathematics, computer science, and cryptography. For the history of Leonardo
sequences, see [A001595] in the On-Line Encyclopedia of Integer Sequences
\cite{OEIS}, and for the properties of Leonardo numbers, we refer to \cite%
{Dijkstra, Catarino, Alp, Shannon, K}.

A relation between Leonardo $p$-numbers and Fibonacci $p$-numbers was given
by%
\begin{equation}
\mathcal{L}_{p,n}=\left( p+1\right) F_{p,n+1}-p  \label{leo1}
\end{equation}%
and a relation between Leonardo $p$-numbers, Fibonacci $p$-numbers and Lucas
$p$-numbers was given by%
\begin{equation}
\mathcal{L}_{p,n}=L_{p,n+p+1}-F_{p,n+p+1}-p.  \label{leo2}
\end{equation}

Recently, Zhong et al.\cite{zhong} have studied a companion sequence of
Leonardo $p$-sequence, called Lucas-Leonardo $p$-sequence, and defined it by
the following non-homogenous relation:%
\begin{equation}
\mathcal{R}_{p,n}=\mathcal{R}_{p,n-1}+\mathcal{R}_{p,n-p-1}+p,\mbox{   }n>p
\label{luc}
\end{equation}%
with initial values $\mathcal{R}_{p,0}=p^{2}+p+1,\mathcal{R}_{p,1}=\cdots =%
\mathcal{R}_{p,p}=1.$ The non-homogenous relation of Lucas-Leonardo $p$%
-sequence can be converted to the following homogenous recurrence relation
\begin{equation}
\mathcal{R}_{p,n}=\mathcal{R}_{p,n-1}+\mathcal{R}_{p,n-p}-\mathcal{R}%
_{p,n-2p-1,}\mbox{   }n>2p.  \label{luc0}
\end{equation}%
A relation between Lucas-Leonardo $p$-numbers and Lucas $p$-numbers is given
by%
\begin{equation}
\mathcal{R}_{p,n}=\left( p+1\right) L_{p,n}-p  \label{luc1}
\end{equation}%
and a relation between Lucas-Leonardo $p$-numbers and Leonardo $p$-numbers
is given by%
\begin{equation}
\mathcal{R}_{p,n}=\left( p+1\right) \mathcal{L}_{p,n}-p\mathcal{L}_{p,n-1}.
\label{luc2}
\end{equation}%
In particular, for $p=1$, it reduces to the Lucas-Leonardo sequence $%
\left
\{ \mathcal{R}_{n}\right \} _{n=0}^{\infty }$ in [A022319]. A
relation between Lucas-Leonardo and Lucas numbers is
\begin{equation}
\mathcal{R}_{n}=2L_{n}-1.  \label{*}
\end{equation}

Another companion sequence of Leonardo sequence, namely the Francois
sequence $\left \{ \mathcal{F}_{n}\right \} _{n=0}^{\infty }$, see [A022318]
in \cite{OEIS}, is studied by Diskaya and Menken \cite{diskaya}. The
Francois numbers also satisfy the same recurrence relation as Lucas-Leonardo
numbers but begins with initial values $\mathcal{F}_{0}=2,\mathcal{F}_{1}=1.$
A relation between Francois numbers, Fibonacci and Lucas numbers is given by%
\begin{equation}
\mathcal{F}_{n}=L_{n}+F_{n+1}-1.  \label{**}
\end{equation}

On the other hand, quaternion algebras finds applications in various fields
such as mathematics, physics, computer science, engineering, and robotics.
It is particularly useful in representing and manipulating spatial rotations
and orientations in three-dimensional space due to its compact and efficient
representation of such transformations. It is also known that the quaternion
algebras are central simple algebras and also clifford algebras. Let $F$ be
a field with characteristic not $2.$ The generalized quaternion algebra over
a field $F$ is defined as:%
\begin{equation}
\begin{array}{c}
\mathbb{Q}_{F}\left( a,b\right) =\left\{ x=x_{1}+x_{2}i+x_{3}j+x_{4}k\mid
x_{1},x_{2},x_{3},x_{4}\in F,\right.  \\
\left. i^{2}=a,j^{2}=b,ij=-ji=k\right\}
\end{array}
\label{q}
\end{equation}%
where $a,b$ are nonzero invertible elements of field $F$. It is clear to see
that when $F=%
\mathbb{R}
$ and $a,b=-1,$ we get the real quaternion algebra. For $x\in \mathbb{Q}%
_{F}\left( a,b\right) ,$ the norm of $x,$ denoted as $N(x)$ and defined as $%
N(x)=x_{1}^{2}-ax_{2}^{2}-bx_{3}^{2}+abx_{4}^{2}.$ We recall that a
generalized quaternion algebra is a division algebra if and only if a
quaternion with a norm of zero is necessarily the zero quaternion.
Otherwise, the algebra is called a split algebra. It is known that real
quaternion algebra is a division algebra and the quaternion algebra over
finite field $%
\mathbb{Z}
_{q},$denoted as $Q_{%
\mathbb{Z}
_{q}}\left( -1,-1\right) $ is a split algebra, where $q$ is an odd prime
integer, see \cite{grau}. Many research studies have concentrated on
quaternions whose components stem from distinctive integer sequences such as
Fibonacci, Lucas, Leonardo sequences, among various others. In particular,
Horadam \cite{horadam} introduced the Fibonacci quaternions by taking
coefficients as Fibonacci numbers. In \cite{savin}, Savin studied special
Fibonacci quaternions in quaternion algebras over finite fields. Recently,
in \cite{tan-diana}, the authors studied Leonardo quaternions over finite
fields. For more information related to quaternion algebras and
Fibonacci-like quaternion sequences we refer to \cite{flaut1, flaut2,
savin1, halici, tan-leung2, dagdeviren} and the references therein.

In this paper, we would like to contribute on this topic by introducing a
new class of quaternions, called Lucas-Leonardo $p$-quaternions. We derive
several fundamental properties of these numbers including recurrence
relations, generating function, and summation formulas. The main motivation
of \ our paper is to consider companion sequences of Leonardo quaternions
and derive some special Lucas-Leonardo and Francois quaternions in
quaternion algebras over finite fields. In particular, we determine the
Lucas-Leonardo quaternions and Francois quaternions which are zero divisors
and invertible elements in the quaternion algebra over certain finite fields.

\section{Lucas-Leonardo $p$-quaternions}

In this section, we introduce a new class of quaternions, namely
Lucas-Leonardo $p$-quaternions, and explore the fundamental properties of
these sequences. Throughout this section, we consider the real quaternion
algebra $\mathbb{Q}_{%
\mathbb{R}
}\left( -1,-1\right) .$ Also, for the simplicity, we denote $I$ as $1+i+j+k.$

\begin{definition}
\label{def1}The $n^{\text{th}}$ Lucas-Leonardo $p$-quaternion is defined by%
\begin{equation*}
\mathbb{Q}\mathcal{R}_{p,n}=\mathcal{R}_{p,n}+\mathcal{R}_{p,n+1}i+\mathcal{R%
}_{p,n+2}j+\mathcal{R}_{p,n+3}k
\end{equation*}%
where $\mathcal{R}_{p,n}$ is the $n^{\text{th}}$ Lucas-Leonardo $p$-number.
\end{definition}

From the definition of Lucas-Leonardo $p$-numbers, we have%
\begin{eqnarray*}
\mathbb{Q}\mathcal{R}_{p,0} &=&\mathcal{R}_{p,0}+\mathcal{R}_{p,1}i+\mathcal{%
R}_{p,2}j+\mathcal{R}_{p,3}k=I+p\left( p+1\right) \\
\mathbb{Q}\mathcal{R}_{p,1} &=&\mathcal{R}_{p,1}+\mathcal{R}_{p,2}i+\mathcal{%
R}_{p,3}j+\mathcal{R}_{p,4}k=I \\
&&\vdots \\
\mathbb{Q}\mathcal{R}_{p,p-3} &=&\mathcal{R}_{p,p-3}+\mathcal{R}_{p,p-2}i+%
\mathcal{R}_{p,p-1}j+\mathcal{R}_{p,p}k=I \\
\mathbb{Q}\mathcal{R}_{p,p-2} &=&\mathcal{R}_{p,p-2}+\mathcal{R}_{p,p-1}i+%
\mathcal{R}_{p,p}j+\mathcal{R}_{p,p+1}k=I+\left( p+1\right) ^{2}k \\
\mathbb{Q}\mathcal{R}_{p,p-1} &=&\mathcal{R}_{p,p-1}+\mathcal{R}_{p,p}i+%
\mathcal{R}_{p,p+1}j+\mathcal{R}_{p,p+2}k \\
&=&1+i+\left( p^{2}+2p+2\right) j+\left( p^{2}+3p+3\right) k \\
&=&I+\left( p+1\right) \left( \left( p+1\right) j+\left( p+2\right) k\right)
\\
\mathbb{Q}\mathcal{R}_{p,p} &=&\mathcal{R}_{p,p}+\mathcal{R}_{p,p+1}i+%
\mathcal{R}_{p,p+2}j+\mathcal{R}_{p,p+3}k \\
&=&1+\left( p^{2}+2p+2\right) i+\left( p^{2}+3p+3\right) j+\left(
p^{2}+5p+5\right) k \\
&=&I+\left( p+1\right) \left( \left( p+1\right) i+\left( p+2\right) j+\left(
p+4\right) k\right) .
\end{eqnarray*}

It is clear to see that from the relation (\ref{luc}) and the definition of
Lucas-Leonardo $p$-quaternions, we have%
\begin{equation}
\mathbb{Q}\mathcal{R}_{p,n}=\mathbb{Q}\mathcal{R}_{p,n-1}+\mathbb{Q}\mathcal{%
R}_{p,n-p-1}+pI,\mbox{   }n>p.  \label{r1}
\end{equation}

For $n>2p,$ from the relation (\ref{luc0}) and the definition of
Lucas-Leonardo $p$-numbers, we have

\begin{equation}
\mathbb{Q}\mathcal{R}_{p,n}=\mathbb{Q}\mathcal{R}_{p,n-1}+\mathbb{Q}\mathcal{%
R}_{p,n-p}-\mathbb{Q}\mathcal{R}_{p,n-2p-1}.  \label{r2}
\end{equation}

Now we state the generating function of the Lucas-Leonardo $p$-quaternions.

\begin{proposition}
The generating function for Lucas-Leonardo $p$-quaternions is%
\begin{equation*}
G\left( x\right) =\frac{I+p\left( p+1\right) +\sum \limits_{n=1}^{p}\left(
\mathbb{Q}\mathcal{R}_{p,n}-\mathbb{Q}\mathcal{R}_{p,n-1}\right) x^{n}+pI%
\frac{x^{p+1}}{1-x}}{1-x-x^{p+1}}.
\end{equation*}
\end{proposition}

\begin{proof}
We consider the formal power series representation for the generating
function of Lucas-Leonardo $p$-quaternions, $G\left( x\right)
=\sum_{n=0}^{\infty }\mathbb{Q}\mathcal{R}_{p,n}x^{n}.$ By using the
relation (\ref{r1}), we get
\begin{equation*}
\left( 1-x-x^{p+1}\right) G\left( x\right) =\sum_{n=0}^{\infty }\mathbb{Q}%
\mathcal{R}_{p,n}x^{n}-\sum_{n=0}^{\infty }\mathbb{Q}\mathcal{R}%
_{p,n}x^{n+1}-\sum_{n=0}^{\infty }\mathbb{Q}\mathcal{R}_{p,n}x^{n+p+1}
\end{equation*}%
\begin{equation*}
=\sum_{n=0}^{p}\mathbb{Q}\mathcal{R}_{p,n}x^{n}-\sum_{n=1}^{p}\mathbb{Q}%
\mathcal{R}_{p,n-1}x^{n}+\sum_{n=p+1}^{\infty }\left( \mathbb{Q}\mathcal{R}%
_{p,n}-\mathbb{Q}\mathcal{R}_{p,n-1}-\mathbb{Q}\mathcal{R}_{p,n-p-1}\right)
x^{n}\text{ \ \ \ \ \ }
\end{equation*}%
\begin{equation*}
=\sum_{n=0}^{p}\mathbb{Q}\mathcal{R}_{p,n}x^{n}-\sum_{n=1}^{p}\mathbb{Q}%
\mathcal{R}_{p,n-1}x^{n}+pI\sum_{n=p+1}^{\infty }x^{n}\text{ \ \ \ \ \ \ \ \
\ \ \ \ \ \ \ \ \ \ \ \ \ \ \ \ \ \ \ \ \ \ \ \ \ \ \ \ \ \ \ \ \ \ \ \ \ \
\ \ \ \ \ \ \ }
\end{equation*}%
\begin{equation*}
=\mathbb{Q}\mathcal{R}_{p,0}+\sum_{n=1}^{p}\left( \mathbb{Q}\mathcal{R}%
_{p,n}-\mathbb{Q}\mathcal{R}_{p,n-1}\right) x^{n}+pI\frac{x^{p+1}}{1-x}\text{
\ \ \ \ \ \ \ \ \ \ \ \ \ \ \ \ \ \ \ \ \ \ \ \ \ \ \ \ \ \ \ \ \ \ \ \ \ \
\ \ \ \ \ \ \ \ \ }
\end{equation*}%
\begin{equation*}
=I+p\left( p+1\right) +\sum_{n=1}^{p}\left( \mathbb{Q}\mathcal{R}_{p,n}-%
\mathbb{Q}\mathcal{R}_{p,n-1}\right) x^{n}+pI\frac{x^{p+1}}{1-x}\text{ \ \ \
\ \ \ \ \ \ \ \ \ \ \ \ \ \ \ \ \ \ \ \ \ \ \ \ \ \ \ \ \ \ \ \ \ \ \ \ }
\end{equation*}%
which gives the desired result.
\end{proof}

Note that from initial values of Lucas-Leonardo $p$-quaternions, we have%
\begin{eqnarray*}
&&\sum \limits_{n=1}^{p}\left( \mathbb{Q}\mathcal{R}_{p,n}-\mathbb{Q}%
\mathcal{R}_{p,n-1}\right) x^{n} \\
&=&\left( \mathbb{Q}\mathcal{R}_{p,1}-\mathbb{Q}\mathcal{R}_{p,0}\right)
x+\left( \mathbb{Q}\mathcal{R}_{p,2}-\mathbb{Q}\mathcal{R}_{p,1}\right)
x^{2}+\cdots +\left( \mathbb{Q}\mathcal{R}_{p,p-3}-\mathbb{Q}\mathcal{R}%
_{p,p-4}\right) x^{p-3} \\
&&+\left( \mathbb{Q}\mathcal{R}_{p,p-2}-\mathbb{Q}\mathcal{R}_{p,p-3}\right)
x^{p-2}+\left( \mathbb{Q}\mathcal{R}_{p,p-1}-\mathbb{Q}\mathcal{R}%
_{p,p-2}\right) x^{p-1}+\left( \mathbb{Q}\mathcal{R}_{p,p}-\mathbb{Q}%
\mathcal{R}_{p,p-1}\right) x^{p} \\
&=&-p\left( p+1\right) x+\left( p+1\right) ^{2}kx^{p-2}+\left( p+1\right)
^{2}\left( j+k\right) x^{p-1}+\left( p+1\right) \left( \left( p+1\right)
i+j+2k\right) x^{p} \\
&=&\left( p+1\right) \left( -px+\left( p+1\right) kx^{p-2}+\left( p+1\right)
\left( j+k\right) x^{p-1}+\left( \left( p+1\right) i+j+2k\right)
x^{p}\right) .
\end{eqnarray*}

Now we provide some relations between Lucas-Leonardo $p$-quaternions,
Leonardo $p$-quaternions, Fibonacci $p$-quaternions, and Lucas $p$%
-quaternions. It's worth noting that Fibonacci $p$-quaternions and Lucas $p$%
-quaternions can be defined similarly to Lucas-Leonardo $p$-quaternions as:
\begin{equation*}
\mathbb{Q}F_{p,n}=F_{p,n}+F_{p,n+1}i+F_{p,n+2}j+F_{p,n+3}k
\end{equation*}%
and%
\begin{equation*}
\mathbb{Q}L_{p,n}=L_{p,n}+L_{p,n+1}i+L_{p,n+2}j+L_{p,n+3}k,
\end{equation*}%
respectively.

\begin{proposition}
\label{p}For $n\geq p,$ we have the following relations:

$%
\begin{array}{c}
\left( i\right) \text{ }\mathbb{Q}\mathcal{R}_{p,n}=\left( p+1\right)
\mathbb{Q}L_{p,n}-pI,%
\end{array}%
$

$%
\begin{array}{c}
\left( ii\right) \text{ }\mathbb{Q}\mathcal{R}_{p,n}=\left( p+1\right)
\mathbb{Q}\mathcal{L}_{p,n}-p\mathbb{Q}\mathcal{L}_{p,n-1},%
\end{array}%
$

$%
\begin{array}{c}
\left( iii\right) \text{ }\mathbb{Q}\mathcal{R}_{p,n}=\mathbb{Q}\mathcal{L}%
_{p,n}+p\mathbb{Q}\mathcal{L}_{p,n-p-1}-p^{2}I,%
\end{array}%
$

$%
\begin{array}{c}
\left( iv\right) \text{ }\mathbb{Q}\mathcal{R}_{p,n}=p\left( \mathbb{Q}%
L_{p,n}-\mathbb{Q}F_{p,n}\right) +\mathbb{Q}L_{p,n+p+1}-\mathbb{Q}%
F_{p,n+p+1}-p\left( p+1\right) I.%
\end{array}%
$
\end{proposition}

\begin{proof}
$\left( i\right) $ From (\ref{luc1}), we have%
\begin{eqnarray*}
\mathbb{Q}\mathcal{R}_{p,n} &=&\mathcal{R}_{p,n}+\mathcal{R}_{p,n+1}i+%
\mathcal{R}_{p,n+2}j+\mathcal{R}_{p,n+3}k \\
&=&\left( \left( p+1\right) L_{p,n}-p\right) +\left( \left( p+1\right)
L_{p,n+1}-p\right) i \\
&&+\left( \left( p+1\right) L_{p,n+2}-p\right) j+\left( \left( p+1\right)
L_{p,n+3}-p\right) k \\
&=&\left( p+1\right) \left( L_{p,n}+L_{p,n+1}i+L_{p,n+2}j+L_{p,n+3}k\right)
-p\left( 1+i+j+k\right) \\
&=&\left( p+1\right) \mathbb{Q}L_{p,n}-pI.
\end{eqnarray*}

$\left( ii\right) $ From (\ref{luc2}), we have%
\begin{eqnarray*}
\mathbb{Q}\mathcal{R}_{p,n} &=&\mathcal{R}_{p,n}+\mathcal{R}_{p,n+1}i+%
\mathcal{R}_{p,n+2}j+\mathcal{R}_{p,n+3}k \\
&=&\left( p+1\right) \mathcal{L}_{p,n}-p\mathcal{L}_{p,n-1}+\left( \left(
p+1\right) \mathcal{L}_{p,n+1}-p\mathcal{L}_{p,n}\right) i \\
&&+\left( \left( p+1\right) \mathcal{L}_{p,n+2}-p\mathcal{L}_{p,n+1}\right)
j+\left( \left( p+1\right) \mathcal{L}_{p,n+3}-p\mathcal{L}_{p,n+2}\right) k
\\
&=&\left( p+1\right) \left( \mathcal{L}_{p,n}+\mathcal{L}_{p,n+1}i+\mathcal{L%
}_{p,n+2}j+\mathcal{L}_{p,n+3}k\right) \\
&&-p\left( \mathcal{L}_{p,n-1}+\mathcal{L}_{p,n}i+\mathcal{L}_{p,n+1}j+%
\mathcal{L}_{p,n+2}k\right) \\
&=&\left( p+1\right) \mathbb{Q}\mathcal{L}_{p,n}-p\mathbb{Q}\mathcal{L}%
_{p,n-1}.
\end{eqnarray*}%
$\left( iii\right) $ By using the relation $\left( ii\right) $ and (\ref{leo}%
), we can easily obtained the second relation.%
\begin{eqnarray*}
\mathbb{Q}\mathcal{R}_{p,n} &=&\left( p+1\right) \mathbb{Q}\mathcal{L}%
_{p,n}-p\mathbb{Q}\mathcal{L}_{p,n-1} \\
&=&\left( p+1\right) \left( \mathcal{L}_{p,n}+\mathcal{L}_{p,n+1}i+\mathcal{L%
}_{p,n+2}j+\mathcal{L}_{p,n+3}k\right) \\
&&-p\left( \mathcal{L}_{p,n-1}+\mathcal{L}_{p,n}i+\mathcal{L}_{p,n+1}j+%
\mathcal{L}_{p,n+2}k\right) \\
&=&p\left( \mathcal{L}_{p,n}-\mathcal{L}_{p,n-1}\right) +p\left( \mathcal{L}%
_{p,n+1}-\mathcal{L}_{p,n}\right) i+p\left( \mathcal{L}_{p,n+2}-\mathcal{L}%
_{p,n+1}\right) j+p\left( \mathcal{L}_{p,n+3}+\mathcal{L}_{p,n+2}\right) k \\
&&+\left( \mathcal{L}_{p,n}+\mathcal{L}_{p,n+1}i+\mathcal{L}_{p,n+2}j+%
\mathcal{L}_{p,n+3}k\right) \\
&=&p\left( \mathcal{L}_{p,n-p-1}-p\right) +p\left( \mathcal{L}%
_{p,n-p}-p\right) i+p\left( \mathcal{L}_{p,n-p+1}-p\right) j+p\left(
\mathcal{L}_{p,n-p+2}-p\right) k \\
&&+\left( \mathcal{L}_{p,n}+\mathcal{L}_{p,n+1}i+\mathcal{L}_{p,n+2}j+%
\mathcal{L}_{p,n+3}k\right) \\
&=&p\left( \mathcal{L}_{p,n-p-1}+\mathcal{L}_{p,n-p}i+\mathcal{L}_{p,n-p+1}j+%
\mathcal{L}_{p,n-p+2}k\right) -p^{2}\left( 1+i+j+k\right) +\mathbb{Q}%
\mathcal{L}_{p,n} \\
&=&p\mathbb{Q}\mathcal{L}_{p,n-p-1}+\mathbb{Q}\mathcal{L}_{p,n}-p^{2}I.
\end{eqnarray*}%
$\left( iv\right) $ By using the relation $\left( iii\right) $ and (\ref%
{leo2}), we obtain the desired result.
\begin{eqnarray*}
\mathbb{Q}\mathcal{R}_{p,n} &=&p\mathbb{Q}\mathcal{L}_{p,n-p-1}+\mathbb{Q}%
\mathcal{L}_{p,n}-p^{2}I. \\
&=&p\left( \mathcal{L}_{p,n-p-1}+\mathcal{L}_{p,n-p}i+\mathcal{L}_{p,n-p+1}j+%
\mathcal{L}_{p,n-p+2}k\right) \\
&&+\left( \mathcal{L}_{p,n}+\mathcal{L}_{p,n+1}i+\mathcal{L}_{p,n+2}j+%
\mathcal{L}_{p,n+3}k\right) -p^{2}I \\
&=&p\left( L_{p,n}-F_{p,n}-p\right) +p\left( L_{p,n+1}-F_{p,n+1}-p\right) i
\\
&&+p\left( L_{p,n+2}-F_{p,n+2}-p\right) j+p\left( \left(
L_{p,n+3}-F_{p,n+3}-p\right) \right) k \\
&&+\left( L_{p,n+p+1}+L_{p,n+p+2}i+L_{p,n+p+3}j+L_{p,n+p+4}k\right) \\
&&-\left( F_{p,n+p+1}+F_{p,n+p+2}i+F_{p,n+p+3}j+F_{p,n+p+4}k\right)
-pI-p^{2}I \\
&=&p\left( L_{p,n}+L_{p,n+1}i+L_{p,n+2}j+L_{p,n+3}k\right) -p\left(
F_{p,n}+F_{p,n+1}i+F_{p,n+2}j+F_{p,n+3}k\right) \\
&&+\mathbb{Q}L_{p,n+p+1}-\mathbb{Q}F_{p,n+p+1}-pI-p^{2}I \\
&=&p\left( \mathbb{Q}L_{p,n}-\mathbb{Q}F_{p,n}\right) +\mathbb{Q}L_{p,n+p+1}-%
\mathbb{Q}F_{p,n+p+1}-p\left( p+1\right) I.
\end{eqnarray*}
\end{proof}

\begin{proposition}
For $n\geq p$, we have%
\begin{equation*}
\sum \limits_{k=0}^{n}\mathbb{Q}\mathcal{R}_{p,k}=\left( p+1\right) \left(
\mathbb{Q}L_{p,n+p+1}-\mathbb{Q}L_{p,p}\right) -p\left( n+1\right) I.
\end{equation*}
\end{proposition}

\begin{proof}
By using the relation $\left( i\right) $ in Proposition \ref{p} and the
summation formula for Lucas $p$-quaternions $\sum \limits_{r=0}^{n}\mathbb{Q}%
L_{p,r}=\mathbb{Q}L_{p,n+p+1}-\mathbb{Q}L_{p,p}$, we get%
\begin{eqnarray*}
\sum \limits_{r=0}^{n}\mathbb{Q}\mathcal{R}_{p,r} &=&\left( p+1\right) \sum
\limits_{r=0}^{n}\mathbb{Q}L_{p,r}-\sum \limits_{r=0}^{n}pI \\
&=&\left( p+1\right) \sum \limits_{r=0}^{n}\mathbb{Q}L_{p,r}-p\left(
n+1\right) I \\
&=&\left( p+1\right) \left( \mathbb{Q}L_{p,n+p+1}-\mathbb{Q}L_{p,p}\right)
-p\left( n+1\right) I.
\end{eqnarray*}
\end{proof}

Next we give a convoluted relation for Lucas-Leonardo $p$-quaternions. To do
this, we need the following relations in \cite[Theorem 3]{abbad} and \cite[%
Property 6]{tuglu}, respectively:%
\begin{equation}
\sum \limits_{t=1}^{p}L_{p,n-t}F_{p,t}=L_{p,n+p}-F_{p,p+1}L_{p,n},
\label{abbad}
\end{equation}%
\begin{equation}
\sum \limits_{t=0}^{n}F_{p,t}=F_{p,n+p+1}-F_{p,p}.\text{ \ \ \ \ \ \ \ \ \ \
\ \ }  \label{tuglu}
\end{equation}

\begin{proposition}
For $n\geq p$, we have%
\begin{equation*}
\sum \limits_{t=1}^{p}\mathbb{Q}\mathcal{R}_{p,n-t}F_{p,t}=\left( p+1\right)
\left( \mathbb{Q}L_{p,n+p}-\mathbb{Q}L_{p,n}\right) -pI\left(
F_{p,2p+1}-1\right) .
\end{equation*}
\end{proposition}

\begin{proof}
By using the relation $\left( i\right) $ in Proposition \ref{p} and (\ref%
{tuglu}), we have%
\begin{equation*}
\sum \limits_{t=1}^{p}\mathbb{Q}\mathcal{R}_{p,n-t}F_{p,t}=\sum%
\limits_{t=1}^{p}\left( \left( p+1\right) \mathbb{Q}L_{p,n-t}-pI\right)
F_{p,t}\text{ \ \ \ \ \ \ \ \ \ \ \ \ \ \ \ \ \ \ \ \ \ \ \ \ \ \ \ \ \ \ \
\ \ \ \ \ \ \ \ \ \ \ \ \ \ \ \ \ \ \ \ \ \ \ \ \ \ \ \ \ \ \ \ \ \ \ \ \ \
\ \ \ \ \ \ }
\end{equation*}%
\begin{equation*}
=\left( p+1\right) \sum \limits_{t=1}^{p}\mathbb{Q}L_{p,n-t}F_{p,t}-pI\sum%
\limits_{t=1}^{p}F_{p,t}\text{ \ \ \ \ \ \ \ \ \ \ \ \ \ \ \ \ \ \ \ \ \ \ \
\ \ \ \ \ \ \ \ \ \ \ \ \ \ \ \ \ }
\end{equation*}%
\begin{equation*}
=\left( p+1\right) \sum \limits_{t=1}^{p}\left(
L_{p,n-t}+L_{p,n-t+1}i+L_{p,n-t+2}j+L_{p,n-t+3}k\right) F_{p,t}-pI\left(
F_{p,2p+1}-F_{p,p}\right) \text{ }
\end{equation*}%
\ \ \ \ \ \ \ \ \ \ \ \ \ \ \ \ \ \ \ \ \ \ \ \ \ \ \ \ \ \ \ \ \ \ \ \ \ \
\ \ \ \ \ \ \ \ \ \ \ \ \ \ \ \ \
\begin{eqnarray*}
&=&\left( p+1\right) \left( \sum \limits_{t=1}^{p}L_{p,n-t}F_{p,t}+\sum
\limits_{t=1}^{p}L_{p,n-t+1}F_{p,t}i\right. \\
&&\left. +\sum
\limits_{t=1}^{p}L_{p,n-t+2}F_{p,t}j+\sum%
\limits_{t=1}^{p}L_{p,n-t+3}F_{p,t}k\right) -pI\left( F_{p,2p+1}-1\right)
\text{ }.\text{ \ \ \ \ \ \ \ \ \ \ }
\end{eqnarray*}%
From the relation (\ref{abbad}) and since $F_{p,p+1}=1$, we get%
\begin{eqnarray*}
\sum \limits_{t=1}^{p}\mathbb{Q}\mathcal{R}_{p,n-t}F_{p,t} &=&\left(
p+1\right) \left( L_{p,n+p}-F_{p,p+1}L_{p,n}\ +\left(
L_{p,n+p+1}-F_{p,p+1}L_{p,n+1}\right) i\right. \\
&&\left. +\left( L_{p,n+p+2}-F_{p,p+1}L_{p,n+2}\ \right) j+\left(
L_{p,n+p+3}-F_{p,p+1}L_{p,n+3}\ \right) k\right) \\
&&-pI\left( F_{p,2p+1}-1\right) \text{ } \\
&=&\left( p+1\right) \left(
L_{p,n+p}+L_{p,n+p+1}i+L_{p,n+p+2}j+L_{p,n+p+3}k\right) \\
&&-\left( L_{p,n}+L_{p,n+1}i+L_{p,n+2}\ j+L_{p,n+3}k\right) \\
&&-pI\left( F_{p,2p+1}-1\right) \text{ } \\
&=&\left( p+1\right) \mathbb{Q}L_{p,n+p}-\left( p+1\right) \mathbb{Q}%
L_{p,n}-pI\left( F_{p,2p+1}-1\right) .
\end{eqnarray*}%
\
\end{proof}

\section{Lucas-Leonardo quaternions and Francois quaternions over finite
fields}

In this section, we consider the quaternion algebra $\mathbb{Q}_{%
\mathbb{Z}
q}\left( -1,-1\right) ,$ for simplicity $\mathbb{Q}_{%
\mathbb{Z}
q}$ where $q$ is an odd prime integer. First, we determine the
Lucas-Leonardo quaternions which are zero divisors in the quaternion algebra
$\mathbb{Q}_{%
\mathbb{Z}
q}$ for $q=3,5,$ and $7$. It's important to note that finding zero divisors
and invertible elements in the Lucas-Leonardo $p$-quaternions is more
challenging task than in the Fibonacci quaternions studied by Savin \cite%
{savin}. This is mainly because the norm of Lucas-Leonardo $p$-quaternions
is more complex. So, we focus on the standard Lucas-Leonardo quaternion case
instead. Similar approach will be used for the Francois quaternions.

Let $\mathbb{Q}\mathcal{R}_{n}$ be the $n$th Lucas-Leonardo quaternion
defined as%
\begin{equation*}
\mathbb{Q}\mathcal{R}_{n}=\mathcal{R}_{n}+\mathcal{R}_{n+1}i+\mathcal{R}%
_{n+2}j+\mathcal{R}_{n+3}k
\end{equation*}%
where the basis $\left \{ 1,i,j,k\right \} $ satisfies the multiplication
rules: $i^{2}=j^{2}=k^{2}=ijk=-1.$ By using the relation (\ref{*}) for the
Lucas-Leonardo quaternion and using the relations $%
L_{n}^{2}+L_{n+1}^{2}=5F_{2n+1},$ $F_{n}^{2}+F_{n+1}^{2}=F_{2n+1},$ $%
F_{n}+F_{n+4}=3F_{n+2},$ $L_{n}F_{n+1}=F_{2n+1}+\left( -1\right) ^{n}F_{1},$
$L_{n}+L_{n+2}=5F_{n+1},F_{n}+F_{n+2}=L_{n+1}$ the norm of Lucas-Leonardo
quaternion can be obtained as follows:%
\begin{eqnarray}
N\left( \mathbb{Q}\mathcal{R}_{n}\right) &=&\mathcal{R}_{n}^{2}+\mathcal{R}%
_{n+1}^{2}+\mathcal{R}_{n+2}^{2}+\mathcal{R}_{n+3}^{2}  \notag \\
&=&\left( 2L_{n}-1\right) ^{2}+\left( 2L_{n+1}-1\right) ^{2}+\left(
2L_{n+2}-1\right) ^{2}+\left( 2L_{n+3}-1\right) ^{2}  \notag \\
&=&4\left( L_{n}^{2}+L_{n+1}^{2}+L_{n+2}^{2}+L_{n+3}^{2}\right) -4\left(
L_{n}+L_{n+1}+L_{n+2}+L_{n+3}\right) +4  \notag \\
&=&20\left( F_{2n+1}+F_{2n+5}\right) -4\left( L_{n+2}+L_{n+4}\right) +4
\notag \\
&=&60F_{2n+3}-20F_{n+3}+4=4\left( 15F_{2n+3}-5F_{n+3}+1\right) .
\label{norm}
\end{eqnarray}

\begin{proposition}
A Lucas-Leonardo quaternion $\mathbb{Q}\mathcal{R}_{n}$ is a zero divisor in
quaternion algebra $\mathbb{Q}_{%
\mathbb{Z}
_{3}}$ if and only if $n\equiv 0,2,3$ $\left( \func{mod}8\right) .$
\end{proposition}

\begin{proof}
A Lucas-Leonardo quaternion $\mathbb{Q}\mathcal{R}_{n}$ is a zero divisor in
quaternion algebra $\mathbb{Q}_{%
\mathbb{Z}
_{3}}$ if and only if $N\left( \mathbb{Q}\mathcal{R}_{n}\right) =\overline{0}
$ in $%
\mathbb{Z}
_{3}.$ By using the relation (\ref{norm}), we have%
\begin{equation*}
F_{n+3}+1\equiv 0\left( \func{mod}3\right) \Leftrightarrow F_{n+3}\equiv
2\left( \func{mod}3\right)
\end{equation*}%
Since the cycle of the Fibonacci numbers is modulo $3$ is%
\begin{equation}
0,1,1,2,0,2,2,1,0,1,1,2,0,2,2,1,\ldots ,  \label{cyclefib}
\end{equation}%
the cycle length of Fibonacci numbers modulo $3$ is $8$. So%
\begin{equation*}
F_{n+3}\equiv 2\left( \func{mod}3\right) \Leftrightarrow n+3\equiv
3,5,6\left( \func{mod}8\right) \Leftrightarrow n\equiv 0,2,3\left( \func{mod}%
8\right) .
\end{equation*}
\end{proof}

\begin{proposition}
All Lucas-Leonardo quaternions $\mathbb{Q}\mathcal{R}_{n}$ are invertible in
quaternion algebra $\mathbb{Q}_{%
\mathbb{Z}
_{5}}$.
\end{proposition}

\begin{proof}
Since $N\left( \mathbb{Q}\mathcal{R}_{n}\right) =\overline{4}\neq \overline{0%
}$ in $%
\mathbb{Z}
_{5},$ all Lucas-Leonardo quaternions $\mathbb{Q}\mathcal{R}_{n}$ are
invertible in quaternion algebra $\mathbb{Q}_{%
\mathbb{Z}
_{5}}.$
\end{proof}

From previous Proposition, there are no Lucas-Leonardo quaternions which are
zero divisors in quaternion algebra $\mathbb{Q}_{%
\mathbb{Z}
_{5}}.$

\begin{proposition}
A Lucas-Leonardo quaternion $\mathbb{Q}\mathcal{R}_{n}$ is a zero divisor in
quaternion algebra $\mathbb{Q}_{%
\mathbb{Z}
_{7}}$ if and only if $n\equiv 0,6,7,9\left( \func{mod}16\right) .$
\end{proposition}

\begin{proof}
A Lucas-Leonardo quaternion $\mathbb{Q}\mathcal{R}_{n}$ is a zero divisor in
quaternion algebra $\mathbb{Q}_{%
\mathbb{Z}
_{7}}$ if and only if $N\left( \mathbb{Q}\mathcal{R}_{n}\right) =\overline{0}
$ in $%
\mathbb{Z}
_{7}.$ By using the relation (\ref{norm}), we have%
\begin{eqnarray*}
N\left( \mathbb{Q}\mathcal{R}_{n}\right) &=&\overline{0}\Leftrightarrow
4\left( 15F_{2n+3}-5F_{n+3}+1\right) \equiv 0\left( \func{mod}7\right) \\
&\Leftrightarrow &4\left( F_{2n+3}+2F_{n+3}+1\right) \equiv 0\left( \func{mod%
}7\right) \\
&\Leftrightarrow &F_{2n+3}+2F_{n+3}+1\equiv 0\left( \func{mod}7\right) \\
&\Leftrightarrow &F_{n+1}^{2}+F_{n+2}^{2}+2\left( F_{n+1}+F_{n+2}\right)
+1\equiv 0\left( \func{mod}7\right) \\
&\Leftrightarrow &\left( F_{n+1}+1\right) ^{2}+\left( F_{n+2}+1\right)
^{2}\equiv 1\left( \func{mod}7\right) .
\end{eqnarray*}%
Since the cycle of Fibonacci numbers modulo $7$ is%
\begin{equation*}
0,1,1,2,3,5,1,6,0,6,6,5,4,2,6,1,\ldots ,
\end{equation*}%
the cycle length of Fibonacci numbers modulo $7$ is $16$. See \cite[A105870]%
{OEIS}. So, we get%
\begin{equation*}
\left( F_{n+1}+1\right) ^{2}+\left( F_{n+2}+1\right) ^{2}\equiv 1\left(
\func{mod}7\right) \Leftrightarrow n\equiv 0,6,7,9\left( \func{mod}16\right)
.
\end{equation*}
\end{proof}

Next, we consider the Francois quaternions over certain finite fields. Let $%
\mathbb{Q}\mathcal{F}_{n}$ be the $n$th Francois quaternion defined as%
\begin{equation*}
\mathbb{Q}\mathcal{F}_{n}=\mathcal{F}_{n}+\mathcal{F}_{n+1}i+\mathcal{F}%
_{n+2}j+\mathcal{F}_{n+3}k
\end{equation*}%
where the basis $\left \{ 1,i,j,k\right \} $ satisfies the real quaternion
multiplication rules. By using the relation (\ref{**}) for the Francois
quaternion, the norm of Francois quaternion can be obtained as follows:%
\begin{eqnarray*}
N\left( \mathbb{Q}\mathcal{F}_{n}\right) &=&\mathcal{F}_{n}^{2}+\mathcal{F}%
_{n+1}^{2}+\mathcal{F}_{n+2}^{2}+\mathcal{F}_{n+3}^{2} \\
&=&\left( L_{n}+F_{n+1}-1\right) ^{2}+\left( L_{n+1}+F_{n+2}-1\right)
^{2}+\left( L_{n+2}+F_{n+3}-1\right) ^{2}+\left( L_{n+3}+F_{n+4}-1\right)
^{2} \\
&=&\left( L_{n}^{2}+L_{n+1}^{2}+L_{n+2}^{2}+L_{n+3}^{2}\right) +\left(
F_{n+1}^{2}+F_{n+2}^{2}+F_{n+3}^{2}+F_{n+4}^{2}\right) +4 \\
&&+2\left( L_{n}F_{n+1}+L_{n+1}F_{n+2}+L_{n+2}F_{n+3}+L_{n+3}F_{n+4}\right)
-2\left( L_{n}+L_{n+1}+L_{n+2}+L_{n+3}\right) \\
&&-2\left( F_{n+1}+F_{n+2}+F_{n+3}+F_{n+4}\right) \\
&=&5F_{2n+1}+5F_{2n+5}+3F_{2n+5}+4 \\
&&+2\left( F_{2n+1}+\left( -1\right) ^{n}+F_{2n+3}+\left( -1\right)
^{n+1}+F_{2n+5}+\left( -1\right) ^{n+2}+F_{2n+7}+\left( -1\right)
^{n+3}\right) \\
&&-2\left( L_{n+2}+L_{n+4}\right) -2\left( F_{n+3}+F_{n+5}\right) \\
&=&5F_{2n+1}+5F_{2n+5}+3F_{2n+5}+4+2\left(
F_{2n+1}+F_{2n+3}+F_{2n+5}+F_{2n+7}\right) \\
&&-2\left( 5F_{n+3}\right) -2\left( F_{n+3}+F_{n+5}\right) \\
&=&7F_{2n+1}+10F_{2n+5}+4+2\left( F_{2n+1}+F_{2n+2}\right) +2\left(
F_{2n+5}+F_{2n+6}\right) -12F_{n+3}-2F_{n+5} \\
&=&9F_{2n+1}+12F_{2n+5}+4+2F_{2n+2}+2\left( F_{2n+4}+F_{2n+5}\right)
-34F_{n+1}-18F_{n} \\
&=&39F_{2n+1}+48F_{2n+2}-34F_{n+1}-18F_{n}+4.
\end{eqnarray*}

\begin{proposition}
A Francois quaternion $Q\mathcal{F}_{n}$ is a zero divisor in quaternion
algebra $Q_{%
\mathbb{Z}
_{3}}$ if and only if $n\equiv 0,1,6$ $\left( \func{mod}8\right) .$
\end{proposition}

\begin{proof}
A Francois quaternion $Q\mathcal{F}_{n}$ is a zero divisor in quaternion
algebra $Q_{%
\mathbb{Z}
_{3}}$ if and only if $N\left( Q\mathcal{F}_{n}\right) =\overline{0}$ in $%
\mathbb{Z}
_{3}.$ From (\ref{cyclefib}), the cycle length of Fibonacci numbers modulo $%
3 $ is $8$. Thus, we have%
\begin{eqnarray*}
N\left( \mathbb{Q}\mathcal{F}_{n}\right) &=&\overline{0}\left( \func{mod}%
3\right) \Leftrightarrow -F_{n+1}+1\equiv 0\left( \func{mod}3\right)
\Leftrightarrow F_{n+1}\equiv 1\left( \func{mod}3\right) \\
&\Leftrightarrow &n+1\equiv 1,2,7\left( \func{mod}8\right) \Leftrightarrow
n\equiv 0,1,6\left( \func{mod}8\right) .
\end{eqnarray*}
\end{proof}

\begin{proposition}
A Francois quaternion $\mathbb{Q}\mathcal{F}_{n}$ is a zero divisor in
quaternion algebra $\mathbb{Q}_{%
\mathbb{Z}
_{5}}$ if and only if $n\equiv 5,8,10,19$ $\left( \func{mod}20\right) .$
\end{proposition}

\begin{proof}
A Francois quaternion $\mathbb{Q}\mathcal{F}_{n}$ is a zero divisor in
quaternion algebra $\mathbb{Q}_{%
\mathbb{Z}
_{5}}$ if and only if $N\left( \mathbb{Q}\mathcal{F}_{n}\right) =\overline{0}
$ in $%
\mathbb{Z}
_{5}.$ By using the norm relation, we have%
\begin{eqnarray}
N\left( \mathbb{Q}\mathcal{F}_{n}\right) &=&\overline{0}\Leftrightarrow
4F_{2n+1}+3F_{2n+2}+F_{n+1}+2F_{n}+4\equiv 0\left( \func{mod}5\right)  \notag
\\
&\Leftrightarrow &-F_{2n+1}-2F_{2n+2}+F_{n+2}+F_{n}+4\equiv 0\left( \func{mod%
}5\right)  \notag \\
&\Leftrightarrow &-F_{2n+3}-F_{2n+2}+F_{n+2}+F_{n}+4\equiv 0\left( \func{mod}%
5\right)  \notag \\
&\Leftrightarrow &F_{2n+3}+F_{2n+2}-F_{n+2}-F_{n}+1\equiv 0\left( \func{mod}%
5\right)  \notag \\
&\Leftrightarrow &F_{2n+4}-F_{n+2}-F_{n}+1\equiv 0\left( \func{mod}5\right) .
\label{cong}
\end{eqnarray}%
We recall that the cycle of Fibonacci numbers modulo $5$ is%
\begin{equation*}
0,1,1,2,3,0,3,3,1,4,0,4,4,3,2,0,2,2,4,1.
\end{equation*}%
So the cycle length of Fibonacci numbers modulo $5$ is $20$. See \cite[%
A082116]{OEIS}. To find $n$ such that the congruence (\ref{cong}) is
satisfied, we need to consider the following cases:

\textbf{Case 1:} If $n\equiv 0\left( \func{mod}20\right) ,$ then $n+2\equiv
2\left( \func{mod}20\right) ,2n+4\equiv 4\left( \func{mod}20\right)
,F_{2n+4}\equiv 3\left( \func{mod}5\right) ,F_{n+2}\equiv 1\left( \func{mod}%
5\right) ,F_{n}\equiv 0\left( \func{mod}5\right) .$ It results $%
F_{2n+4}-F_{n+2}-F_{n}+1\equiv 3\left( \func{mod}5\right) ,$therefore the
congruence (\ref{cong}) is not satisfied.

\textbf{Case 2:} If $n\equiv 1\left( \func{mod}20\right) ,$ then $n+2\equiv
3\left( \func{mod}20\right) ,2n+4\equiv 6\left( \func{mod}20\right)
,F_{2n+4}\equiv 3\left( \func{mod}5\right) ,F_{n+2}\equiv 2\left( \func{mod}%
5\right) ,F_{n}\equiv 1\left( \func{mod}5\right) .$ It results $%
F_{2n+4}-F_{n+2}-F_{n}+1\equiv 1\left( \func{mod}5\right) ,$ therefore the
congruence (\ref{cong}) is not satisfied.

\textbf{Case 3: } If $n\equiv 2\left( \func{mod}20\right) ,$ then $n+2\equiv
4\left( \func{mod}20\right) ,2n+4\equiv 8\left( \func{mod}20\right)
,F_{2n+4}\equiv 1\left( \func{mod}5\right) ,F_{n+2}\equiv 3\left( \func{mod}%
5\right) ,F_{n}\equiv 1\left( \func{mod}5\right) .$ It results $%
F_{2n+4}-F_{n+2}-F_{n}+1\equiv 3\left( \func{mod}5\right) ,$ therefore the
congruence (\ref{cong}) is not satisfied.

\textbf{Case 4: } If $n\equiv 3\left( \func{mod}20\right) ,$ then $n+2\equiv
5\left( \func{mod}20\right) ,2n+4\equiv 10\left( \func{mod}20\right)
,F_{2n+4}\equiv 0\left( \func{mod}5\right) ,F_{n+2}\equiv 0\left( \func{mod}%
5\right) ,F_{n}\equiv 2\left( \func{mod}5\right) .$ It results $%
F_{2n+4}-F_{n+2}-F_{n}+1\equiv 4\left( \func{mod}5\right) ,$ therefore the
congruence (\ref{cong}) is not satisfied.

\textbf{Case 5: } If $n\equiv 4\left( \func{mod}20\right) ,$ then $n+2\equiv
6\left( \func{mod}20\right) ,2n+4\equiv 12\left( \func{mod}20\right)
,F_{2n+4}\equiv 4\left( \func{mod}5\right) ,F_{n+2}\equiv 3\left( \func{mod}%
5\right) ,F_{n}\equiv 3\left( \func{mod}5\right) .$ It results $%
F_{2n+4}-F_{n+2}-F_{n}+1\equiv 4\left( \func{mod}5\right) ,$ therefore the
congruence (\ref{cong}) is not satisfied.

\textbf{Case 6: } If $n\equiv 5\left( \func{mod}20\right) ,$ then $n+2\equiv
7\left( \func{mod}20\right) ,2n+4\equiv 14\left( \func{mod}20\right)
,F_{2n+4}\equiv 2\left( \func{mod}5\right) ,F_{n+2}\equiv 3\left( \func{mod}%
5\right) ,F_{n}\equiv 0\left( \func{mod}5\right) .$ So, the congruence (\ref%
{cong}) is satisfied. Therefore, we have $N\left( \mathbb{Q}\mathcal{F}%
_{n}\right) =\overline{0}\Leftrightarrow n\equiv 5\left( \func{mod}20\right)
.$

\textbf{Case 7: } If $n\equiv 6\left( \func{mod}20\right) ,$ then $n+2\equiv
8\left( \func{mod}20\right) ,2n+4\equiv 16\left( \func{mod}20\right)
,F_{2n+4}\equiv 2\left( \func{mod}5\right) ,F_{n+2}\equiv 1\left( \func{mod}%
5\right) ,F_{n}\equiv 3\left( \func{mod}5\right) .$ It results $%
F_{2n+4}-F_{n+2}-F_{n}+1\equiv 4\left( \func{mod}5\right) ,$ therefore the
congruence (\ref{cong}) is not satisfied.

\textbf{Case 8: } If $n\equiv 7\left( \func{mod}20\right) ,$then $n+2\equiv
9\left( \func{mod}20\right) ,2n+4\equiv 18\left( \func{mod}20\right)
,F_{2n+4}\equiv 4\left( \func{mod}5\right) ,F_{n+2}\equiv 4\left( \func{mod}%
5\right) ,F_{n}\equiv 3\left( \func{mod}5\right) .$ It results $%
F_{2n+4}-F_{n+2}-F_{n}+1\equiv 3\left( \func{mod}5\right) ,$ therefore the
congruence (\ref{cong}) is not satisfied.

\textbf{Case 9: } If $n\equiv 8\left( \func{mod}20\right) ,$ then $n+2\equiv
10\left( \func{mod}20\right) ,2n+4\equiv 0\left( \func{mod}20\right)
,F_{2n+4}\equiv 0\left( \func{mod}5\right) ,F_{n+2}\equiv 0\left( \func{mod}%
5\right) ,F_{n}\equiv 1\left( \func{mod}5\right) .$ So, the congruence (\ref%
{cong}) is satisfied. Therefore, we have $N\left( \mathbb{Q}\mathcal{F}%
_{n}\right) =\overline{0}\Leftrightarrow n\equiv 8\left( \func{mod}20\right)
.$

\textbf{Case 10: } If $n\equiv 9\left( \func{mod}20\right) ,$ then $%
n+2\equiv 11\left( \func{mod}20\right) ,2n+4\equiv 2\left( \func{mod}%
20\right) ,F_{2n+4}\equiv 1\left( \func{mod}5\right) ,F_{n+2}\equiv 4\left(
\func{mod}5\right) ,F_{n}\equiv 4\left( \func{mod}5\right) .$ It results $%
F_{2n+4}-F_{n+2}-F_{n}+1\equiv 4\left( \func{mod}5\right) ,$ therefore the
congruence (\ref{cong}) is not satisfied.

\textbf{Case 11: } If $n\equiv 10\left( \func{mod}20\right) ,$ then $%
n+2\equiv 12\left( \func{mod}20\right) ,2n+4\equiv 4\left( \func{mod}%
20\right) ,F_{2n+4}\equiv 3\left( \func{mod}5\right) ,F_{n+2}\equiv 4\left(
\func{mod}5\right) ,F_{n}\equiv 0\left( \func{mod}5\right) .$ So, the
congruence (\ref{cong}) is satisfied. Therefore, we have $N\left( Q\mathcal{F%
}_{n}\right) =\overline{0}\Leftrightarrow n\equiv 10\left( \func{mod}%
20\right) .$

\textbf{Case 12: } If $n\equiv 11\left( \func{mod}20\right) ,$ then $%
n+2\equiv 13\left( \func{mod}20\right) ,2n+4\equiv 6\left( \func{mod}%
20\right) ,F_{2n+4}\equiv 3\left( \func{mod}5\right) ,F_{n+2}\equiv 3\left(
\func{mod}5\right) ,F_{n}\equiv 4\left( \func{mod}5\right) .$ It results $%
F_{2n+4}-F_{n+2}-F_{n}+1\equiv 2\left( \func{mod}5\right) ,$ therefore the
congruence (\ref{cong}) is not satisfied.

\textbf{Case 13: } If $n\equiv 12\left( \func{mod}20\right) ,$ then $%
n+2\equiv 14\left( \func{mod}20\right) ,2n+4\equiv 8\left( \func{mod}%
20\right) ,F_{2n+4}\equiv 1\left( \func{mod}5\right) ,F_{n+2}\equiv 2\left(
\func{mod}5\right) ,F_{n}\equiv 4\left( \func{mod}5\right) .$ It results $%
F_{2n+4}-F_{n+2}-F_{n}+1\equiv 1\left( \func{mod}5\right) ,$ therefore the
congruence (\ref{cong}) is not satisfied.

\textbf{Case 14: } If $n\equiv 13\left( \func{mod}20\right) ,$ then $%
n+2\equiv 15\left( \func{mod}20\right) ,2n+4\equiv 10\left( \func{mod}%
20\right) ,F_{2n+4}\equiv 0\left( \func{mod}5\right) ,F_{n+2}\equiv 0\left(
\func{mod}5\right) ,F_{n}\equiv 3\left( \func{mod}5\right) .$ It results $%
F_{2n+4}-F_{n+2}-F_{n}+1\equiv 3\left( \func{mod}5\right) ,$ therefore the
congruence (\ref{cong}) is not satisfied.

\textbf{Case 15: } If $n\equiv 14\left( \func{mod}20\right) ,$then $%
n+2\equiv 16\left( \func{mod}20\right) ,2n+4\equiv 12\left( \func{mod}%
20\right) ,F_{2n+4}\equiv 4\left( \func{mod}5\right) ,F_{n+2}\equiv 2\left(
\func{mod}5\right) ,F_{n}\equiv 2\left( \func{mod}5\right) .$ It results $%
F_{2n+4}-F_{n+2}-F_{n}+1\equiv 1\left( \func{mod}5\right) ,$ therefore the
congruence (\ref{cong}) is not satisfied.

\textbf{Case 16: } If $n\equiv 15\left( \func{mod}20\right) ,$ then $%
n+2\equiv 17\left( \func{mod}20\right) ,2n+4\equiv 14\left( \func{mod}%
20\right) ,F_{2n+4}\equiv 2\left( \func{mod}5\right) ,F_{n+2}\equiv 2\left(
\func{mod}5\right) ,F_{n}\equiv 0\left( \func{mod}5\right) .$ It results $%
F_{2n+4}-F_{n+2}-F_{n}+1\equiv 1\left( \func{mod}5\right) ,$ therefore the
congruence (\ref{cong}) is not satisfied.

\textbf{Case 17: } If $n\equiv 16\left( \func{mod}20\right) ,$ then $%
n+2\equiv 18\left( \func{mod}20\right) ,2n+4\equiv 16\left( \func{mod}%
20\right) ,F_{2n+4}\equiv 2\left( \func{mod}5\right) ,F_{n+2}\equiv 4\left(
\func{mod}5\right) ,F_{n}\equiv 2\left( \func{mod}5\right) .$ It results $%
F_{2n+4}-F_{n+2}-F_{n}+1\equiv 2\left( \func{mod}5\right) ,$ therefore the
congruence (\ref{cong}) is not satisfied.

\textbf{Case 18: } If $n\equiv 17\left( \func{mod}20\right) ,$ then $%
n+2\equiv 19\left( \func{mod}20\right) ,2n+4\equiv 18\left( \func{mod}%
20\right) ,F_{2n+4}\equiv 4\left( \func{mod}5\right) ,F_{n+2}\equiv 1\left(
\func{mod}5\right) ,F_{n}\equiv 2\left( \func{mod}5\right) .$ It results $%
F_{2n+4}-F_{n+2}-F_{n}+1\equiv 2\left( \func{mod}5\right) ,$ therefore the
congruence (\ref{cong}) is not satisfied.

\textbf{Case 19: } If $n\equiv 18\left( \func{mod}20\right) ,$ then $%
n+2\equiv 0\left( \func{mod}20\right) ,2n+4\equiv 0\left( \func{mod}%
20\right) ,F_{2n+4}\equiv 0\left( \func{mod}5\right) ,F_{n+2}\equiv 0\left(
\func{mod}5\right) ,F_{n}\equiv 4\left( \func{mod}5\right) .$ It results $%
F_{2n+4}-F_{n+2}-F_{n}+1\equiv 2\left( \func{mod}5\right) ,$ therefore the
congruence (\ref{cong}) is not satisfied.

\textbf{Case 20: } If $n\equiv 19\left( \func{mod}20\right) ,$ then $%
n+2\equiv 1\left( \func{mod}20\right) ,2n+4\equiv 2\left( \func{mod}%
20\right) ,F_{2n+4}\equiv 1\left( \func{mod}5\right) ,F_{n+2}\equiv 1\left(
\func{mod}5\right) ,F_{n}\equiv 1\left( \func{mod}5\right) .$ So, the
congruence (\ref{cong}) is satisfied. Therefore, we have $N\left( \mathbb{Q}%
\mathcal{F}_{n}\right) =\overline{0}\Leftrightarrow n\equiv 19\left( \func{%
mod}20\right) .$

Thus we obtain the desired result.
\end{proof}

\section{Conclusion}

Quaternions find broad applications across numerous fields, notably in
physics and mathematics. Therefore, it is valuable to explore and analyze
the properties of specific types of quaternions. Numerous studies have
focused on quaternions whose components originate from special integer
sequences, including Fibonacci, Lucas, Leonardo sequences, among others.
Here, we introduce a new class of quaternions which can be seen as a
companion quaternion sequence of Leonardo $p$-quaternions. We investigate
basic properties of these quaternions and give some applications over finite
fields. The results of this paper can be summarized as follows:

\begin{itemize}
\item We derive a new class of quaternions which takes Lucas-Leonardo $p$%
-numbers as a components. When $p=1,$ the Lucas-Leonardo $p$-quaternions
reduce to the Lucas-Leonardo quaternions. We obtain generating function,
recurrence relations and sum identities for these new quaternion class.

\item We obtain that

\begin{itemize}
\item a Lucas-Leonardo quaternion $\mathbb{Q}\mathcal{R}_{n}$ is a zero
divisor in quaternion algebra $\mathbb{Q}_{%
\mathbb{Z}
_{3}}$ if and only if $n\equiv 0,2,3$ $\left( \func{mod}8\right) ,$

\item all Lucas-Leonardo quaternions are invertible in quaternion algebra $%
\mathbb{Q}_{%
\mathbb{Z}
_{5}},$

\item a Lucas-Leonardo quaternion $\mathbb{Q}\mathcal{R}_{n}$ is a zero
divisor in quaternion algebra $\mathbb{Q}_{%
\mathbb{Z}
_{7}}$ if and only if $n\equiv 0,6,7,9\left( \func{mod}16\right) .$
\end{itemize}

\item We obtain that

\begin{itemize}
\item A Francois quaternion $Q\mathcal{F}_{n}$ is a zero divisor in
quaternion algebra $Q_{%
\mathbb{Z}
_{3}}$ if and only if $n\equiv 0,1,6$ $\left( \func{mod}8\right) ,$

\item A Francois quaternion $\mathbb{Q}\mathcal{F}_{n}$ is a zero divisor in
quaternion algebra $\mathbb{Q}_{%
\mathbb{Z}
_{5}}$ if and only if $n\equiv 5,8,10,19$ $\left( \func{mod}20\right) .$
\end{itemize}
\end{itemize}

\end{document}